\documentclass{amsart}
\usepackage{amssymb}
\usepackage{mathrsfs}
\usepackage{stmaryrd} 
\usepackage{bm}
\usepackage{bbm}
\usepackage{oldgerm}
\usepackage[francais]{babel}
\usepackage{amscd}
\usepackage[T1]{fontenc}
\usepackage[latin1]{inputenc}
\usepackage[all]{xy}
\usepackage[pagebackref=true]{hyperref}
\renewcommand*{\backref}[1]{}
\renewcommand*{\backrefalt}[4]{[{\tiny%
    \ifcase #1 Non cit\'e.%
          \or Cit\'e \`a la page~#2.%
          \else Cit\'e aux pages #2.%
    \fi%
    }]}

\usepackage[final]{pdfpages} 
\usepackage{upref}
\usepackage{setspace}
\usepackage{xcolor}

\hypersetup{
    colorlinks,
    linkcolor={black},
    citecolor={black},
    urlcolor={black}
}

\newtheorem{theo}[subsection]{Th\'{e}or\`{e}me}
\newtheorem{prop}[subsection]{Proposition}
\newtheorem{cor}[subsection]{Corollaire}

\theoremstyle{definition}

\newtheorem{defi}[subsection]{D\'{e}finition}

\numberwithin{equation}{subsection}

\newcommand{\cD}{{\mathscr D}}
\newcommand{\cE}{{\mathscr E}}
\newcommand{\cF}{{\mathscr F}}

\newcommand{\cI}{{\mathscr I}}

\newcommand{\cK}{{\mathscr K}}

\newcommand{\cP}{{\mathscr P}}
\newcommand{\co}{{\mathscr O}}

\newcommand{\cT}{{\mathscr T}}
\newcommand{\cH}{{\mathscr H}}

\newcommand{\lbr}{\begin{bmatrix}}
\newcommand{\rbr}{\end{bmatrix}}

\begin{document}

\title[$q$-déformation de la théorie de Hodge non-abélienne en caractéristique positive]
{Sur une  $q$-déformation locale de la théorie de Hodge non-abélienne en caractéristique positive}
 \author{ Michel Gros}
 
\address{M.G. CNRS UMR 6625, IRMAR, Universit\'{e} de Rennes 1,
Campus de Beaulieu, 35042 Rennes cedex, France}

\email{michel.gros@univ-rennes1.fr}

\begin{abstract}

Pour $p$ un nombre premier et $q$ une racine $p$-ième   non triviale de 1, nous présentons les principales étapes de la construction d'une $q$-déformation locale de la ``correspondance de Simpson  en caractéristique $p$'' dégagée par Ogus et Vologodsky  en 2005. La construction est basée sur  l'équivalence de Morita entre un anneau d'opérateurs différentiels $q$-déformés et   son centre. Nous expliquons aussi les liens espérés entre cette construction et celles introduites récemment par Bhatt et Scholze. Pour alléger l'exposition, nous nous limitons au cas de la dimension 1.\\

For $p$ a prime number and  $q$  a non trivial $p$th root  of   1,  we present the main steps of the construction of  a local  $q$-deformation of the ``Simpson correspondence  in   characteristic $p$'' found by Ogus and Vologodsky in 2005. The construction is based on the Morita-equivalence between a ring of  $q$-twisted differential operators and its center. We also explain the expected relations between this construction and those recently done by Bhatt and Scholze. For the sake of readability, we limit ourselves to the case of dimension 1.

\end{abstract}

\maketitle

\setcounter{tocdepth}{1}
\tableofcontents

\section{Introduction}
\subsection{}\label{intro1}

Ogus et Vologodsky ont dégagé dans \cite{ov} un   analogue en caractéristique $p>0$ de la théorie de Hodge non-abélienne, {\emph{i.e.}}  de la correspondance de Simpson complexe. Soient  $\tilde{S}$ un schéma plat   sur ${\Bbb{Z}}/p^{2}$,  $\tilde{X}$, $\tilde{X'}$ deux $\tilde{S}$-schémas lisses de réduction modulo $p$ notées $X$ et $X'$,    et  ${\rm{{\tilde{F}}}} : {\tilde{X}} \rightarrow {\tilde{X'}}$ un $\tilde{S}$-morphisme. Supposons que ces données constituent un  relèvement au-dessus de $\tilde{S}$ du morphisme de Frobenius relatif ${\rm{F}}_{X/S} : X \rightarrow X'$ associé à $X$ vu comme schéma au-dessus de $S= \tilde{S}\times_{{\Bbb{Z}}/p^{2}} {\Bbb{Z}}/p$. Elles permettent alors à Ogus et Vologodsky d'étendre (\cite{ov}, Thm. 2.8, Thm. 2.26) aux $\co_X$-modules munis d'une connexion intégrable dont la $p$-courbure est supposée seulement {\emph{quasi}}-nilpotente  à la fois le théorème de descente de Cartier  (\cite{katz}, Thm. 5.1)  et  l'existence d'une  décomposition du complexe de De Rham obtenu  par Deligne et Illusie (\cite{di}, Rem. 2.2(ii)) induisant l'opération de Cartier  (\cite{katz}, Thm. 7.2). L'exposé  oral d'A. Abbes et le nôtre ont été consacrés aux travaux d'Oyama (\cite{oyama}), Shiho (\cite{shiho}) et Xu   (\cite{xu}) qui ont permis de relever ``modulo $p^n$''   cette correspondance d'Ogus et Vologodsky. C'est ici   une  autre direction qui est explorée.


\subsection{}\label{intro2}
Sans rapport avec  ce qui précède, Bhatt, Morrow et Scholze   ont  
dégagé  (\cite{bms}, Thm. 1.8) un raffinement entier des théorèmes standards de comparaison entre cohomologies cristalline,  de De Rham et étale $p$-adique pour un schéma  formel  propre  et lisse  sur l'anneau des entiers d'une extension non-archimédienne algébriquement close de ${\Bbb{C}}_{p}$. Dans l'élaboration de celui-ci apparait un relèvement de l'isomorphisme de Cartier (\cite{bms}, Thm. 8.3) sur la cohomologie d'un objet (\cite{bms}, Def. 8.1) d'une certaine catégorie dérivée.  Dans des situations géométriques locales bien adaptées (\cite{bms}, 8.5) auxquelles les auteurs se ramènent pour  établir cet isomorphisme, l'existence de ce relèvement  découle de l'étude    de  certains  $q$-{\emph{complexes de  De Rham}} (\cite{bms}, 7.7) avec  $q$  une racine $p$-ième non triviale de l'unité dans ${\Bbb{C}}_{p}$. Ces derniers sont de vrais complexes  qui ``réalisent'' (\cite{bms}, \S8) les objets des catégories dérivées évoquées ci-dessus. Ils ont  eux-aussi une cohomologie se calculant  par un relèvement  de l'opération de Cartier (\cite{scholze}, Prop. 3.4, (iii) ; voir aussi \cite{prid}, Prop. 2.8) qui   explique donc localement  l'existence de la précédente. Il nous semble plausible que l'extension espérée du théorème de comparaison entier (\cite{bms}, Thm. 1.8) à des coefficients non constants  (\cite{tsuji}) donne quelque intérêt à essayer d'expliciter une {\emph{ $q$-déformation locale de la théorie de Hodge non-abélienne}  incluant  l'étude de ce type de complexes et les propriétés de leurs cohomologies. Cela devrait    peut-être éclaircir un peu d'éventuels  liens entre les  théories  \cite{ov} et \cite{bms}  puis,  ultérieurement,  ceux avec la correspondance de Simpson $p$-adique (\cite{agt}).

\subsection{}\label{intro3}
Le but de ce rapport est d'esquisser, dans ces situations géométriques locales  bien adaptées, une telle variante. Les deux résultas principaux sont, d'une part une $q$-déformation locale de la correspondance développée par Ogus-Vologodsky (\ref{sci}) et, d'autre part,  sa compatibilité aux  cohomologies naturelles du but et de la source (\ref{qcomp}).   Un corollaire facile (\ref{qcart}) de notre résultat est l'existence de l'opération de Cartier ``relevée''  (\cite{scholze}, Prop. 3.4 ; voir aussi \cite{prid}, 2.2). Dans le type de situation géométrique  que nous considérons,   les résultats principaux d'Ogus-Vologodsky qui nous intéressent ici découlent immédiatement d'une équivalence de Morita, à savoir celle associée à la neutralisation d'une algèbre d'opérateurs différentiels vue comme algèbre d'Azumaya sur son centre (\cite{ov}, Thm. 2.11). Nous   développons  simplement  un $q$-analogue de tout le tableau. Plusieurs choix doivent être faits lors des constructions et il  est peu probable que ces résultats puissent se globaliser par les techniques standards,  purement schématiques, de recollement   (voir par exemple \cite{scholze}, Conj. 1.1  et {\emph{infra}}, \cite{prid}, 2.2 et 3.4 pour  une discussion de problèmes analogues, \cite{bms}, Rem. 8.4).

\subsection{}\label{intro3'}
Pour pallier ces difficultés de globalisation des $q$-complexes de De Rham, Bhatt et Scholze ont introduit très récemment dans \cite{bs}, à beaucoup d'autres fins aussi (dont celle  de réinterpréter les théorèmes de comparaison entiers évoqués plus haut ainsi que les décompositions de Hodge-Tate, ...),   de nouvelles techniques, en particulier celles du site prismatique et du site $q$-cristallin. Ils utilisent pour ce faire la théorie des $\delta$-anneaux et leurs avancées  fournissent  pour nous un espoir   de montrer  l'indépendance de tout choix auxiliaire (en particulier d'une coordonnée) dans nos constructions, au moins à isomorphisme près. De toute façon, vu la généralité du cadre dans lequel ils se placent et le  potentiel d'applications, il nous a paru  indispensable d'en tenir compte  et de reconsidérer avec leurs nouveaux outils  les questions que nous nous posions au moment de la conférence et de la première version de cet article  puis d'indiquer les  progrès réalisés depuis lors.


\subsection{}\label{intro4}
Ce rapport ne contient  pas de démonstrations, pour lesquelles on renvoie à \cite{glq} et à \cite{glq2}. Nous insistons plutôt ici sur la mise en parallèle des théories modulo $p$ (\cite{glq}) et $q$-déformées (\cite{glq2}) en spécialisant cette dernière au cadre familier (notations, hypothèses, terminologie, ...) de la théorie de Hodge $p$-adique, ce qui en allège très largement la présentation. Nous résumons très succintement tout d'abord au \S2   les principales étapes suivies dans \cite{glq}  pour établir  la neutralisation ({\emph{loc. cit.}}, Thm. 4.13) mentionnée ci-dessus. La seule nouveauté par rapport à  \cite{glq}  est le résultat de comparaison cohomologique \ref{comp}. Nous passons ensuite au \S3, après avoir précisé le cadre géométrique,   à la définition des opérateurs différentiels $q$-déformés. Bien qu'on puisse les définir plus directement (cf. \ref{op6}), c'est par un processus de dualité et donc via la définition de parties principales $q$-déformées que nous procédons afin de pouvoir  raisonner comme dans la théorie modulo $p$.  Dans le  \S4, nous déterminons le centre de l'algèbre des opérateurs différentiels $q$-déformés et montrons comment on peut diviser l'action induite par le ``Frobenius''  (\ref{frobdiv}) sur les modules de parties principales $q$-déformées.  Que ceci soit possible est pour l'instant  l'aspect le plus miraculeux de toute cette théorie. Nous en déduisons enfin la neutralisation (\ref{frobdiv31}) d'une complétion centrale    de l'algèbre des opérateurs différentiels $q$-déformés. Le \S5 reformule alors l'équivalence de Morita standard qu'on déduit de cette neutralisation en termes de modules munis d'une $q$-dérivation quasi-nilpotente et de modules de Higgs quasi-nilpotents (\ref{sci}) et les conséquences cohomologiques (\ref{qcomp}). Nous terminons enfin au \S6 par quelques observations et questions en relation avec \cite{bs}.

\subsection{}\label{intro6}

Ces résultats sont le fruit d'une collaboration avec B. Le Stum et A. Quir\'os que l'auteur dégage de toute responsabilité pour les erreurs ou imprécisions qui pourraient apparaitre.  L'auteur remercie très sincèrement  la Fondation Simons  et les organisateurs  de la session {\emph{Simons Symposium on $p$-adic Hodge Theory}}  (8-12 Mai 2017), Bhargav Bhatt et Martin Olsson,  de lui avoir donné l'opportunité d'avancer sur toutes les questions soulevées par ce projet.

 \section{Rappels sur la théorie d'Ogus et Vologodsky}\label{1}
 
\subsection{}\label{ov1}
 Nous conservons dans cette section les notations et hypothèses de \ref{intro1} résumées par les deux diagrammes suivants  
 \begin{equation}
\xymatrix{
\tilde{X}\ar[r]^{\rm{{\tilde{F}}}}\ar[rd] & \tilde{X'}\ar[d] &  X\ar[r]^{{\rm{F}}_{X/S}}\ar[rd] &  X'\ar[d]\\
&\tilde{S}\ar[d]&&S\ar[d]\\
&{\Bbb{Z}}/p^{2}&&{\Bbb{Z}}/p }
\end{equation}
dont celui de droite est donc la réduction modulo $p$ de celui de gauche. Dans ce qui suit, nous allégerons la notation ${\rm{F}}_{X/S}$ en simplement  ${\rm{F}}$ mais en attirant bien l'attention du lecteur sur le fait que cet allègement n'est pas tout à fait compatible avec les notations adoptées dans  \cite{glq} (dans {\it{ loc. cit.}}, ${\rm{F}}$ est noté $F_X$  et $F$  y désigne le Frobenius {\em{absolu}} de $X$). Nous supposerons de plus $S$ noethérien pour raccourcir la preuve de \ref{comp}. 

\subsection{}\label{ov1bis}
 Pour alléger les notations, nous noterons ici simplement $
 \cD_{X}$  la $\co_X$-algèbre  $\cD_{X}^{(0)}$ des opérateurs  différentiels de $X/S$ de niveau $m=0$  introduite par  Berthelot \cite{ber}, 2.2.1 et utilisée dans \cite{glq}, Def. 2.5, parfois dénommée {\it{algèbre des  opérateurs PD-différentiels}} ou {\it{algèbre des opérateurs différentiels cristallins}}. Elle est engendrée par $\co_X$ et par les $S$-dérivations de $\co_X$ (cf. \cite{ber}, p. 218, Rem. (i)). Nous noterons  ${\rm{Z}}\cD_{X}
 $ (resp. ${\rm{Z}}\co_X$) le centre de $\cD_{X}
 $  (resp. le  centralisateur dans $\cD_{X}
 $ de sa sous-algèbre $\co_X$). Nous noterons enfin ${\rm{S}}(\cT_{X'})$ la $\co_{X'}$-algèbre (graduée) symétrique du $\co_{X'}$-module  $\cT_{X'}$ des fonctions sur le fibré  cotangent de $X'/S$. L'application de $p$-courbure permet (cf. par exemple \cite{glq}, Prop. 3.6) de construire  un isomorphisme de $\co_X$-algèbres   
\begin{equation}\label{pcourb}
c : {\rm{S}}(\cT_{X'}) \stackrel{\sim} \rightarrow {\rm{F}}_{*}  {\rm{Z}}\cD_{X}
  \,\, ; \,\,D \in \cT_{X'} \mapsto D^p - D^{[p]}.
\end{equation}
 On peut, de même (cf. {\emph{loc. cit.}}), identifier ${\rm{Z}}\co_X$ à ${\rm{F}} ^{*} {\rm{S}}(\cT_{X'}) = \co_{X}\otimes_{\co_{X'}}  {\rm{S}}(\cT_{X'})$.

\subsection{}\label{ov2}
L'algèbre $\cD_{X}$ agit de manière naturelle de façon $\co_{X'}$-linéaire sur $\co_X$. Soit   $\cK_X$ le noyau de la surjection canonique $\cD_{X} \rightarrow \cE nd_{\co_{X'}}(\co_X)$. C'est un idéal bilatère de  $\cD_{X}$. Nous noterons  ${\widehat{\cD_{X}}}$   (resp. ${\widehat{{\rm{Z}}\cD_{X}}}$, resp. ${\widehat{{\rm{S}}(\cT_{X'})}}$, resp. ${\widehat{{\rm{Z}}\co_X}}$, resp. $\co_{X}\otimes_{\co_{X'}} {\widehat{ {\rm{S}}(\cT_{X'})}}$) le  complété   adique de $\cD_{X}$ (resp. ${\rm{Z}}\cD_{X}$,  resp. ${\rm{S}}(\cT_{X'})$, resp. ${\rm{Z}}\co_X$, resp. $\co_{X}\otimes_{\co_{X'}} {\rm{S}}(\cT_{X'})$)  relativement à l'idéal bilatère  $\cK_X$ (resp. $\cK_X \cap {\rm{S}}(\cT_{X'}) $, resp. $\cK_X \cap {\rm{Z}}\co_X$,   resp. $\cK_X \cap (\co_{X}\otimes_{\co_{X'}}  {\rm{S}}(\cT_{X'}))$).

\subsection{}\label{ov3}
Dans cette situation restrictive d'existence de ${\rm{{\tilde{F}}}}$, plusieurs  des résultats généraux de  \cite{ov} (e.g. Thm. 2.8) découlent immédiatement du résultat suivant (\cite{glq}, Thm. 4.13) que nous avons appris  de P. Berthelot et dont nous rappelerons brièvement le principe de preuve ci-dessous (\ref{ov7}, \ref{ov8}).

\begin{theo}\label{sci}
Toute donnée de ($\tilde{X}$, $\tilde{X'}$,  ${\rm{{\tilde{F}}}} : {\tilde{X}} \rightarrow {\tilde{X'}}$) comme précédemment définit canoniquement un isomorphisme de $\co_K$-algèbres
\begin{equation}\label{sci1}
{\widehat{\cD_{X}}} \stackrel{\sim} \rightarrow \cE nd_{{\widehat{{\rm{S}}(\cT_{X'})}}}(\co_{X}\otimes_{\co_{X'}} {\widehat{ {\rm{S}}(\cT_{X'})}}).
\end{equation}
\end{theo}
On remarquera, en prévision  de (\ref{frodiv2}), que le but de (\ref{sci1}) est simplement $\cE nd_{{\widehat{{\rm{Z}}\cD_{X}}}}({\widehat{{\rm{Z}}\co_X}})$ .
\subsection{}\label{ov4}
Un lemme classique d'algèbre linéaire (\cite{glq}, Lem. 5.6) montre alors que les anneaux ${\widehat{\cD_{X}}}$ et ${\widehat{{\rm{S}}(\cT_{X'})}}$ sont, d'une manière complètement explicite,   équivalents  au sens de Morita : les deux foncteurs suivants entre  les catégories de modules sur ces anneaux, ${\bf{Mod}}\,( {\widehat{\cD_{X}}})$  et  ${\bf{Mod}}\,({\widehat{{\rm{S}}(\cT_{X'})}}$, correspondantes sont quasi-inverses l'un de l'autre
\begin{equation}\label{cov1}
{\Bbb{H}} :  {\bf{Mod}}\,( {\widehat{\cD_{X}}})  \rightarrow  {\bf{Mod}}\,({\widehat{{\rm{S}}(\cT_{X'})}}) \,; \,
\cE  \mapsto  \cH om_{{\widehat{\cD_{X}}}}({\rm{F}} ^{*} {\widehat{{\rm{S}}(\cT_{X'})}}, \cE),
\end{equation}
\begin{equation}\label{cov2}
{\Bbb{M}} :  {\bf{Mod}}\, ({\widehat{{\rm{S}}(\cT_{X'})}})   \rightarrow  {\bf{Mod}}\,({\widehat{\cD_{X}}}) \,; \,
\cF\mapsto  \cF \otimes_{{\widehat{{\rm{S}}(\cT_{X'})}}}{\rm{F}} ^{*} {\widehat{{\rm{S}}(\cT_{X'})}}.
\end{equation}

Ce résultat fournit, une fois réinterprété   (cf. \cite{glq}, Prop. 5.2) les objets de ces catégories le résultat suivant
\begin{theo}[]\label{equivaSimpson}
Toute donnée de ($\tilde{X}$, $\tilde{X'}$,  ${\rm{{\tilde{F}}}} : {\tilde{X}} \rightarrow {\tilde{X'}}$) comme précédemment définit canoniquement une équivalence entre la catégorie  des $\co_{X}$-modules munis d'une connexion intégrable de $p$-courbure quasi-nilpotente (cf. \cite{glq}, Prop. 5.5) et   la catégorie des $\co_{X'}$-modules munis d'un  champ de Higgs quasi-nilpotent (cf. \cite{glq}, Prop. 5.4).  
\end{theo}
On vérifie que dans cette équivalence $\co_X$ muni de sa connexion canonique $d$ correspond à $\co_{X'}$ muni du champ de Higgs nul.


\subsection{}\label{ov5}
Le {\emph{complexe de Higgs}} d'un $\co_{X'}$-module de Higgs  $\cF \in  {\bf{Mod}}\, ({\widehat{{\rm{S}}(\cT_{X'})}}) $  est, par définition, le complexe (avec $\cF$ placé en degré 0)
\begin{equation}
0 \rightarrow \cF \stackrel{\theta} \longrightarrow \cF \otimes_{\co_{X'}} \Omega^{1}_{X'} \stackrel{(-)\wedge \theta} \longrightarrow \cF \otimes_{\co_{X'}} \Omega^{2}_{X'} \stackrel{(-)\wedge \theta} \longrightarrow ...
\end{equation}
avec $\theta$ l'application $\co_{X'}$-linéaire  provenant de la structure naturelle de ${\rm{S}}(\cT_{X'})$-module sur $\cF$ et, pour alléger,  $\Omega^{i}_{X'}$  le $\co_S$-module des différentielles {\em{relatives}} de de degré $i$ de $X'/S$ (noté  $\Omega^{i}_{X'/S}$ lorsque une ambigu\"{\i}té est possible). On en donnera ci-dessous (\ref{kosz1}) une autre description.
Il résulte facilement de cette équivalence la
\begin{prop}[]\label{comp}
Si $\cE \in {\bf{Mod}}\,( {\widehat{\cD_{X}}}) $ et $\cF \in    {\bf{Mod}}\, ({\widehat{{\rm{S}}(\cT_{X'})}})$ se correspondent par l'équivalence ci-dessus, alors l'image directe par  ${\rm{F}} $ du complexe de De Rham de $\cE $ est quasi-isomorphe au complexe de Higgs de $\cF$.
\end{prop}
Le principe de démonstration est le suivant. Pour calculer ${\Bbb{R}}\cH om_{{\widehat{\cD_{X}}}}(\co_X, \cE)$, on utilise la résolution     de Spencer de $\co_{X}$ par des $\co_X$-localement libres sur $\cD_{X} $
\begin{equation}\label{spen}
 [...\rightarrow \cD_{X}^{(0)} \otimes_{\co_X} \wedge^{2}\cT_{X}\rightarrow \cD_{X}^{(0)} \otimes_{\co_X} \cT_{X} \rightarrow \cD_{X}]\rightarrow   \co_X \rightarrow 0.
 \end{equation}
On tensorise alors  la partie entre crochets par ${\widehat{\cD_{X}}}$ en préservant l'exactitude de \ref{spen} car ${\widehat{\cD_{X}}}$ est plat sur $\cD_{X} $ puisque c'est le complété 
de $\cD_{X}$ relativement à  un idéal bilatère engendré par une suite centralisante. On a alors, notant $\Omega^{\bullet}_{X}$ pour alléger   le complexe  $\Omega^{\bullet}_{X/S}$ des différentielles relatives de $X/S$, des isomorphismes
\begin{equation}\label{res1}
{\Bbb{R}}\cH om_{{\widehat{\cD_{X}}}}(\co_X, \cE) \simeq \cH om_{{\widehat{\cD_{X}}}}({\widehat{\cD_{X}}} \otimes_{\co_X} \wedge^{\bullet}\cT_{X}, \cE) \simeq \cE \otimes_{\co_X} \Omega^{\bullet}_{X }
\end{equation}
Pour le complexe de Higgs, on utilise  la résolution de Koszul  de $\co_{X'}$
\begin{equation}\label{kosz}
 [... \rightarrow {\rm{S}}(\cT_{X'}) \otimes_{\co_{X'}} \wedge^{2}\cT_{X'} \rightarrow {\rm{S}}(\cT_{X'}) \otimes_{\co_{X'}} \cT_{X'} \rightarrow {\rm{S}}(\cT_{X'})]\rightarrow   \co_{X'} \rightarrow 0.
 \end{equation}
que l'on tensorise par ${\widehat{{\rm{S}}(\cT_{X'})}}$ au-dessus de ${\rm{S}}(\cT_{X'})$ en la laissant exacte. On obtient
\begin{equation}\label{kosz1}
{\Bbb{R}}\cH om_{{\widehat{{\rm{S}}(\cT_{X'})}}}(\co_{X'}, \cF) \simeq \cH om_{{\widehat{{\rm{S}}(\cT_{X'})}}}({\widehat{{\rm{S}}(\cT_{X'})}}\otimes_{\co_{X'}} \wedge^{\bullet}\cT_{X'}, \cF) \simeq \cF \otimes_{\co_{X'}} \Omega^{\bullet}_{X'}. \end{equation}
 Les deux foncteurs dérivés \ref{res1} et \ref{kosz1} pouvant se calculer à l'aide de résolutions injectives du second argument, la proposition s'ensuit grâce à l'équivalence de catégories donnée par ${\Bbb{H}}$ et ${\Bbb{M}}$.

\subsection{}\label{ov6}
Notons ici que, par définition de $\cK_X$ et grâce au lemme d'algèbre linéaire qu'on vient d'évoquer,  les anneaux  $\cD_{X}/\cK_X$ et   $\co_{X'}$ sont   équivalents au sens de Morita: c'est, réinterprété dans ce langage,  le classique théorème de descente de Cartier (\cite{katz}, Thm. 7.2). D'autre part, la proposition \ref{comp} fournit exactement, une fois précisé les isomorphismes,  la décomposition du complexe de De Rham obtenue par Deligne-Illusie (\cite{di}, Rem. 2.2(ii)).
  
\subsection{}\label{ov7}
La démonstration du théorème \ref{sci} procède par dualité. Soient $\cI \subset \co_{X\times_{S}X}$ l'idéal définissant l'immersion diagonale $X\hookrightarrow X\times_{S}X$, $\cP_{X} $ son {\it{enveloppe à puissances divisées}} (notée  $\cP_{X/S,(0)}$  dans \cite{glq}, 2.4), ${\overline{\cI}} \subset \cP_{X  }$ le {\it{PD-idéal}} engendré par $\cI$ et, pour $n$ un entier $\geq 0$,  $\cP_{X  }^{n} = \cP_{X  }/ \cI^{[n+1]}$. On a, par définition, 
\begin{equation} 
\cD_{X,n}   =  \cH om_{\co_{X}}(\cP_{X  }^{n}, \co_X) \,\,;\,\, \cD_{X }  = \cup_{n\geq0} \cD_{X,n} .
\end{equation}
On remarque alors qu'on a simplement un isomorphisme
\begin{equation} 
{\widehat{\cD_{X} }} \stackrel{\sim} \rightarrow  \cH om_{\co_{X}}(\cP_{X  } , \co_X).  
\end{equation}
et que ${\widehat{\cD_{X} }} $ n'est autre que ce qui est classiquement appelé l'{\it{algèbre des opérateurs hyper-PD-différentiels}}.
D'autre part, notons $\Gamma (\Omega^{1}_{X' })$ la $\co_{X'}$-algèbre (graduée) à puissances divisées canoniquement associée au $\co_{X'}$-module $\Omega^{1}_{X }$ (\cite{glq}, Thm. 1.2). Une  vérification d'algèbre linéaire  (\cite{glq}, preuve de Thm. 4.13) fournit un isomorphisme
\begin{equation} 
\cE nd_{{\widehat{{\rm{S}}(\cT_{X'})}}}(\co_{X}\otimes_{\co_{X'}} {\widehat{ {\rm{S}}(\cT_{X'})}}) \stackrel{\sim} \rightarrow  \cH om_{\co_{X}}(\co_{X\times_{X'}X} \otimes_{\co_{X}'} \Gamma (\Omega^{1}_{X' }) , \co_X)  
\end{equation}
de sorte que le théorème \ref{sci} se réduit à la construction  d'un isomorphisme d'algèbres de Hopf
\begin{equation}\label{isofond}
 \co_{X\times_{X'}X} \otimes_{\co_{X'}} \Gamma (\Omega^{1}_{X' })   \stackrel{\sim} \rightarrow \cP_{X }.
 \end{equation}

\subsection{}\label{ov8}
Cette construction procède selon les principales étapes suivantes :

\subsubsection{}\label{ov8step1} 
L'application canonique $\cI \rightarrow \cP_{X  }$ ; $f \rightarrow f^{[p]}$  composée avec la projection canonique $\cP_{X } \rightarrow \cI\cP_{X }$ est une application ${\rm{F}}^{*} $-linéaire nulle sur $\cI^2$ (\cite{glq}, Lem. 3.1). Elle induit donc (\cite{glq}, Prop. 3.2) par passage au quotient et linéarisation une application $\co_X$-linéaire  
\begin{equation} \label{Frobdiv1}
{\rm{F}}^{*} \Omega^{1}_{X'}  \rightarrow \cP_{X }/\cI\cP_{X } 
\end{equation}
\subsubsection{}\label{ov8step2} 
L'application (\ref{Frobdiv1}) s'étend en un isomorphisme de $\co_X$-algèbres à puissances divisées (\cite{glq}, Prop. 3.3)
\begin{equation} \label{Frobdiv2}
{\rm{F}}^{*} \Gamma (\Omega^{1}_{X' })   \stackrel{\sim} \rightarrow \cP_{X }/\cI\cP_{X }.
\end{equation}
\subsubsection{}\label{ov8step3} 
La donnée de ($\tilde{X}$, $\tilde{X'}$,  ${\rm{{\tilde{F}}}} : {\tilde{X}} \rightarrow {\tilde{X'}}$) permet de factoriser le morphisme (\ref{Frobdiv1}) en un morphisme de $\co_X$-modules
\begin{equation} \label{Frobdiv3}
{\rm{F}}^{*} \Omega^{1}_{X' }    \rightarrow \cP_{X }.
\end{equation}
C'est l'application {\emph{Frobenius divisé}},   notée $\frac{1}{p!} {\rm{{\tilde{F}}}}^{*}$ dans (Prop. 4.8, \cite{glq}).
On prendra garde ici que la surjection  canonique $\cP_{X } \rightarrow \cP_{X }/\cI  \cP_{X }$  n'est pas compatible aux puissances divisées.
\subsubsection{}\label{ov8step4} 
L'application (\ref{Frobdiv3}) s'étend en un  morphisme  de $\co_X$-{\emph{algèbres  à puissances divisées}} (\cite{glq}, Prop. 4.8)
\begin{equation} \label{Frobdiv4}
{\rm{F}}^{*} \Gamma (\Omega^{1}_{X' })   \rightarrow \cP_{X }
\end{equation}
factorisant l'isomorphisme \ref{Frobdiv2}.
\subsubsection{}\label{ov8step5} 
L'application (\ref{Frobdiv4}) s'étend canoniquement en un  isomorphisme  de $\co_X$-algèbres  de Hopf (\cite{glq}, Prop. 4.13)
\begin{equation} \label{Frobdiv5}
 \co_{X\times_{X'}X} \otimes_{\co_{X'}}  {\rm{F}}^{*} \Gamma (\Omega^{1}_{X' })   \stackrel{\sim} \rightarrow \cP_{X }. 
\end{equation}
C'est l'isomorphisme (\ref{isofond}) recherché.

\subsection{}\label{ov9}
Il  peut être  utile au lecteur de savoir   que si l'on composait la projection canonique  $\cP_{X }  \rightarrow \cP_{X }/ \cI \cP_{X }$ avec l'inverse de (\ref{Frobdiv2}) et qu'on dualisait l'application obtenue, on retrouverait la composée ${\rm{S}}(\cT_{X'}) \rightarrow  {\rm{F}}_{ *} \cD_{X}$ de l'application de $p$-courbure $c$ (\ref{pcourb}) et de l'application canonique ${\rm{Z}}\cD_{X}  \hookrightarrow \cD_{X} $.

\section{Opérateurs différentiels $q$-déformés}\label{}

\subsection{}\label{op1}
Soient $R$ un anneau commutatif   supposé muni d'un relèvement qu'on notera ici simplement $\rm{F}$ du Frobenius absolu de $R/p$ et $q \in R$. Soient également $A$ une $R$-algèbre   munie d'un morphisme étale $f : R[t] \rightarrow A$ ({\emph{i.e.}} d'un {\emph{framing}} au sens de \cite{scholze}, \S3 ; \cite{bms}, \S8, ...). On munit $R[t]$ des deux morphismes de $R$-algèbres   $\sigma$ et ${\rm{F}}^{*}$ induits par $\sigma(t) = qt$ et ${\rm{F}}^{*}(t)=t^p$. On supposera également, par simplicité, qu'il existe  deux morphismes de $R$-algèbres notés encore $\sigma : A \rightarrow A$ et ${\rm{F}}^{*} :  A' := R _ {\nwarrow\rm{F}} \otimes _{  R}A \rightarrow A$,  tels que  respectivement $\sigma(x) = qx$ et ${\rm{F}}^{*}(1 \otimes x)=x^p$ avec $x:= f(t)$ (élément parfois appelé   {\emph{coordonnée}} sur $A$).  Signalons immédiatement, pour fixer les idées, un exemple particulièrement intéressant pour nous où une telle situation se manifeste. Soient  ${\overline{{\Bbb{Q}}}}_{p}$ une cl\^oture algébrique  de ${\Bbb{Q}}_{p}$,  $q \in {\overline{{\Bbb{Q}}}}_{p} $ une racine $p$-ième de 1 non triviale, $K$ l'extension finie totalement ramifiée de ${\Bbb{Q}}_{p}$ engendrée par $q$ et  $R:=\co_K$ l'anneau des entiers de $K$ muni de $\rm{F} = {\rm{Id}}_{R}$. Alors, la simple donnée de $f$ {\em{étale}}  comme ci-dessus  et des arguments standards suffisent à produire, par passage à la complétion $p$-adique de $A$, une situation comme précédemment pour cette dernière.\\

On s'est limité au cadre de la dimension 1 mais tout ce qui précède et suit vaut en dimension supérieure. On a également fixé une fois pour toutes ce dont on aura besoin mais  les données ne seront utilisés qu'au fur et à mesure (la donnée de Frobenius  n'est pas requise avant \S4).

\subsection{}\label{op2}
Pour $u$ une indéterminée  et $n$ un entier $\geq 0$, on pose $(n)_u= \frac{u^{n}-1}{u-1} \in {\Bbb{Z}}[u]$ ; $(n)_{u}! = \prod_{i=1}^{n} (i)_u   \in {\Bbb{Z}}[u]$ ; $\left(\begin{array}{c}n \\k\end{array}\right)_{u} = \frac{(n)_{u}!}{(k)_{u}!(n-k)_{u}!}  \in {\Bbb{Z}}[u]$.  Si maintenant $q$ est un élément de $R$ comme dans \ref{op1}, les notations  $(n)_q$ ; $(n)_{q}! $ ; $\left(\begin{array}{c}n \\k\end{array}\right)_{q}$ signifient qu'on a évalué les quantités précédentes en $u=q$ afin d'obtenir des éléments de $R$.  Ayant à éviter plus bas une possible confusion avec la notation standard des puissances divisées par des crochets, nous avons adopté la notation $(n)_q$ avec des parenthèses plutôt que la notation $[n]_q$ de (\cite{gl} ou \cite{scholze},  \S1). 
\subsection{}\label{op3}
Soient $A$ comme dans  \ref{op1} et  fixons $y \in A$. On va définir tout d'abord l'analogue $q$-déformé des sections de $\cP_{X}$ (\ref{ov7}) sur un ouvert affine dans ce cadre. Soit $A \langle  \xi   \rangle_{q,y} $   le $A$-module libre de générateurs abstraits notés $\xi^{[n]_{q,y}}$ avec $n\in {\Bbb{N}}$ (cf. \cite{glq2}, \S2). On abrégera, lorsqu'aucune confusion n'en résulte, $\xi^{[0]_{q,y}}$ en 1 et $\xi^{[1]_{q,y}}$ en $\xi$. On notera $I^{[n+1]_{q,y}}$ le sous-$A$-module libre de $A \langle  \xi   \rangle_{q,y} $ engendré par les $\xi^{[k]_{q,y}}$ avec $k>n$. 
\begin{prop}[\cite{glq2}, Prop. 2.2]\label{multdef}
Soient $m,n \in {\Bbb{N}}$. La règle de multiplication
\begin{equation}\label{ }
\xi^{[n]}. \xi^{[m]} = \sum_{i=0}^{{\rm{min}}(m,n)} (-1)^{i} q^{\frac{i(i-1)}{2} }\left(\begin{array}{c}m+n-i \\m\end{array}\right)_{q}\left(\begin{array}{c}m \\i\end{array}\right)_{q} y^{i}\xi^{[m+n-i]}
\end{equation}
permet de munir $A \langle  \xi   \rangle_{q,y} $ une structure de $A$-algèbre commutative et unitaire. L'ensemble $I^{[n+1]}$ est un idéal de $A \langle  \xi   \rangle_{q,y} $.
\end{prop}
On dira alors que $A \langle  \xi   \rangle_{q,y} $ est l'{\emph{anneau des polynômes sur $A$ à puissances divisées $q$-déformées}}. \\

Cette terminologie est justifiée par l'égalité, valide pour tout $n\in {\Bbb{N}}$,   
\begin{equation}
(n)_{q}!\xi^{[n]} = \prod_{i=0}^{n-1}\ (\xi+(i)_{q}y)  =: \xi^{(n)}
\end{equation}
et le fait que les $\xi^{(n)}$ forment, pour $n\in {\Bbb{N}}$, une base de $A[\xi]$ (\cite{glq2}, Lem. 1.1).

\subsection{}\label{op4}
Soit encore  $A$ comme dans  \ref{op1}. Supposons désormais que  $y = (1-q)x \in A$.  Le lecteur remarquera que lorsque $q=1$, l'algèbre $A \langle  \xi   \rangle_{q,y} $ n'est autre que la $A$-algèbre des polynômes à puissances divisées usuelles en $\xi$.
\begin{defi}[\cite{glq2}, Def. 4.2]\label{ }
Soit $n \in {\Bbb{N}}$. Le  $A$-module  des {\em{parties principales $q$-déformées}} de $A$ d'ordre au plus $n$ (et de niveau 0) et le  $A$-module des parties principales $q$-déformées de $A$  sont, respectivement, 
\begin{equation}
{\rm{P}}^{(0)}_{A/R, \sigma, n} =  A \langle  \xi   \rangle_{q,y}/  I^{[n+1]_{q,y}},
\end{equation}
\begin{equation}
{\rm{P}}^{(0)}_{A/R, \sigma} = \underset{\underset{n \in {\Bbb{N}}}{\longleftarrow}}{\lim}\ A \langle  \xi   \rangle_{q,y}/  I^{[n+1]_{q,y}}.
\end{equation}
\end{defi}
Dans la suite, nous allégerons les notations ${\rm{P}}^{(0)}_{A/R, \sigma, n}$ et ${\rm{P}}^{(0)}_{A/R, \sigma}$ en ${\rm{P}} _{A , \sigma, n}$ et ${\rm{P}} _{A , \sigma}$ respectivement.
\subsection{}\label{op5} 
On conserve les hypothèses et notations de \ref{op4}. Rappelons qu'il découle de (\cite{lq}, Prop. 2.10) qu'il existe un unique endomorphisme $R$-linéaire   $\partial_{\sigma}$ de $A$ tel que, pour tous $z_1, z_2 \in A$, on ait
\begin{equation}
\partial_{\sigma}(z_1 z_2) = z_1 \partial_{\sigma}(z_2) +  \sigma(z_1) \partial_{\sigma}(z_2),
\end{equation}
{\emph{i.e.}} une $\sigma$-{\emph{dérivation}} canonique.
Afin de construire les opérateurs différentiels $q$-déformés par dualité, nous aurons besoin de la définition suivante. 
\begin{defi}[\cite{glq2}, Def. 4.5]\label{ }
L'application de Taylor $q$-déformée (de niveau 0) est l'application
\begin{equation}\label{taylor}
  \cT: A  \rightarrow {\rm{P}} _{A , \sigma}
\end{equation}
définie par $\cT (z) = \sum_{k=0}^{+\infty} \partial^{k}_{\sigma}(z)\xi^{[k]}$ pour tout $z \in A$.  
\end{defi}
On peut en fait définir $\cT$ de manière plus formelle (cf. \cite{glq2}, Def. 4.5)  et vérifier que c'est un  morphisme d'anneaux, puis  décrire cette application grâce à $\partial_{\sigma}$ comme on vient de le faire.

Si maintenant $M$ est un $A$-module à gauche, l'écriture ${\rm{P}}_{A, \sigma,n} \otimes '_{R} M$ signifie que nous regardons ${\rm{P}} _{A, \sigma,n}$ comme un $A$-module via l'application $\cT$ (\ref{taylor}). Autrement dit, pour tous $z \in A, s\in M, k \in {\Bbb{N}}$, on a :
\begin{equation}\label{ }
\xi^{[k]} \otimes ' zs = \cT(z)\xi^{[k]} \otimes ' s
\end{equation} 
Ceci permet de définir, pour chaque $n \in {\Bbb{N}}$, 
\begin{equation}\label{opdt}
{\rm{D}}^{(0)}_{A, \sigma,n} = {\rm{Hom}}_{A}( {\rm{P}} _{A, \sigma, n} \otimes '_{A} {A}, A).
\end{equation} 
 
Pour $n \in {\Bbb{N}}$, ces $A$-modules forment un système inductif et permettent donc de considérer
\begin{equation}\label{defop}
{\rm{D}}^{(0)}_{A, \sigma } =  \underset{\underset{n \in {\Bbb{N}}}{\longrightarrow}}{\lim} \,\,  {\rm{D}}^{(0)}_{A, \sigma,n}.
\end{equation} 


On vérifie alors que   la comultiplication
\begin{equation}\label{ }
 {\rm{P}} _{A, \sigma}  \rightarrow  {\rm{P}} _{A, \sigma} \otimes '_{A}  {\rm{P}} _{A, \sigma }
\end{equation}
définie par $\xi^{[n]} \mapsto \sum_{i=0}^{n} \xi^{[n-i]} \otimes ' \xi^{[i]}$ permet de munir ${\rm{D}} _{A, \sigma }^{(0)}$ d'une structure d'anneau (cf. \cite{glq2}, Prop. 5.6).

\subsection{}\label{op6} 
L'anneau (\ref{defop}) ainsi construit par dualité n'est autre (cf. \cite{glq2}, Prop. 5.7) que l'extension de Ore  ${\rm{D}}_{A/R, \sigma }$ de $A$ par $\sigma$ et $\partial_{\sigma}$, c'est-à-dire le $A$-module libre de générateurs abstraits $\partial_{\sigma}^k$ ($k\geq 0$)  avec la règle de commutation $\partial_{\sigma} z = \sigma(z)\partial_{\sigma} + \partial_{\sigma}(z)$ pour tout $z \in A$. Dans la suite, on utilisera cette notation ${\rm{D}}_{A, \sigma }$ pour l'anneau (\ref{defop}) si aucune confusion n'en résulte.

 \section{$p$-courbure et Frobenius divisé $q$-déformés}\label{pcour}
\subsection{}\label{pcour0}
On garde dans tout ce \S \, les notations et hypothèses de \ref{op1} et l'on suppose de plus que $(p)_q=0$ dans $R$ et que $R$ est {\emph{$q$-divisible}}, c'est-à-dire  (cf. \cite{glq2}, 0.3) que pour  tout $ m \in {\Bbb{N}} $, $(m)_q$  est inversible dans $R$ s'il est non nul.  Ces deux conditions (qui ne sont pas nécessaires simultanément dans tous les énoncés) sont réalisés, par exemple, dans  le cas  où  $R= \co_K$ (\ref{op1})  et  $q \neq 1$  une racine $p$-ième de l'unité. On continue de poser  $y=(1-q)x$ comme dans  \ref{op4}.\\

Soient ${\rm{Z}}{\rm{D}}_{A, \sigma }$  le centre de l'anneau ${\rm{D}}_{A, \sigma }$  et ${\rm{Z}}{A}_{A, \sigma }$ le  centralisateur dans ${\rm{D}}_{A, \sigma }$ de sa sous-algèbre $A$.  \\
\begin{prop}[\cite{glq2}, Prop. 6.3, Def 6.5]\label{ }
Il existe une unique application $A$-linéaire de $A$-algèbres
\begin{equation}\label{ }
A[\theta] \rightarrow {\rm{D}}_{A, \sigma } \,\, ; \,\, \theta \mapsto \partial_{\sigma}^{p}
\end{equation}
dite de $p$-courbure  $q$-déformé (ou simplement $p$-courbure tordue). Elle induit un isomorphisme de $A$-algèbres entre $A[\theta]$ et ${\rm{Z}}{A}_{A, \sigma }$  et de $A'$-algèbres entre $A'[\theta]$ et ${\rm{Z}}{\rm{D}}_{A, \sigma }$.
\end{prop}

Cette application est construite par dualité à partir des applications canoniques ${\rm{P}}_{A, \sigma, np} \rightarrow {\rm{P}}_{A, \sigma, np}/(\xi)$.

\subsection{}\label{pcour1}
L'analogue $q$-déformé du calcul local crucial permettant de prouver l'existence de l'isomorphisme  (\ref{Frobdiv2})  est l'énoncé suivant.
\begin{prop}[\cite{glq2}, Def 2.5, Thm. 2.6]\label{ }
L'unique application $A$-linéaire 
\begin{equation}\label{frobeniustordu}
A \langle  \omega   \rangle_{1,y^p} \rightarrow A \langle  \xi   \rangle_{q,y} \,\, ; \,\, \omega^{[k]} \mapsto \xi^{[pk]}
\end{equation}
est appelée  Frobenius divisé $q$-déformé (ou simplement Frobenius divisé tordu). C'est un homomorphisme d'anneaux induisant un isomorphisme de $A$-algèbres
\begin{equation}\label{ }
A \langle  \omega   \rangle_{1,y^p} \stackrel{\sim} \rightarrow A \langle  \xi   \rangle_{q,y}/(\xi).
\end{equation}
\end{prop}

\subsection{}\label{pcour2}
Pour $n$ et $i$ des entiers $\geq 0$, on définit (cf. \cite{glq2}, Def 7.4, Prop. 7.9)  des polynômes $A_{n,i}(u), B_{n,i}(u) \in {\Bbb{Z}}[u]$ par les formules
\begin{equation}\label{ }
A_{n,i}(u) := \sum_{j=0}^{n}(-1)^{n-j}u^{\frac{p(n-j)(n-j-1}{2}}  \left(\begin{array}{c}n \\j\end{array}\right)_{u^p}\left(\begin{array}{c}pj\\i\end{array}\right)_{u}
\end{equation}
et
\begin{equation}\label{defB}
(n)_{u}!A_{n,i}(u)   = (n)_{u^p}! (p)_{u}^{n}B_{n,i}(u)
\end{equation}
Les polynômes $A_{n,i}(u)$ s'introduisent naturellement dans la description de l'action induite par ${\rm{F}}$ sur les modules de parties principales $q$-déformées (cf. \cite{glq2}, Prop. 7.5).  La possibilité de définir les polynômes $B_{n,i}(u)$ vient, quant à elle, de l'examen des coefficients des $A_{n,i}(u)$.
\begin{prop}[\cite{glq2}, Prop. 7.12]\label{frobdiv}
L'application
\begin{equation}\label{frodiv1}
[{\rm{F}}^*] :  A' \langle  \omega   \rangle_{1,y} \rightarrow A \langle  \xi   \rangle_{q,y}\,\, ;\,\, [{\rm{F}}^*](\omega^{[n]}) = \sum_{i=n}^{pn}B_{n,i}(q)x^{pn-i}\xi^{[i]}
\end{equation}
est un homomorphisme d'anneaux. 
\end{prop}
Grâce à cette  application, on montre, comme  pour  \ref{Frobdiv5} : 
\begin{prop}[\cite{glq2}, Prop. 7.13]\label{frodiv2}
L'application $[{\rm{F}} ^*]$ induit un  morphisme de $A$-algèbres
\begin{equation}\label{ }
[{\rm{F}}^*] : A[\xi]/(\xi^{(p)_{q,y}}) \,\, \otimes_{A'} \,\, A' \langle  \omega   \rangle_{1,y}  \stackrel{\sim}  \rightarrow A \langle  \xi   \rangle_{q,y},
\end{equation}
qui est un isomorphisme.  
\end{prop}
C'est l'analogue $q$-déformé du calcul local crucial permettant de prouver l'existence  de l'isomorphisme  (\ref{Frobdiv5}) (\cite{glq}, Thm. 4.13).

\subsection{}\label{pcour4}
Par dualité, on   déduit de \ref{frodiv2} la $q$-déformation  suivante de (\cite{glq}, Prop. 4.8).
\begin{prop}[\cite{glq2}, Prop. 8.1]\label{frodiv }
L'application $[{\rm{F}}^*]$  induit, par dualité, un morphisme de $A$-modules
\begin{equation}\label{frobd}
{\rm{\Phi}}_{A, \sigma} : {\rm{D}}_{A, \sigma } \rightarrow {\rm{Z}}{A}_{A, \sigma } \hookrightarrow  {\rm{D}}_{A, \sigma } \,\,; \,\, \partial_{\sigma}^{n} \mapsto \sum_{k=0}^{n}B_{k,n}(q)x^{pk-n}\partial_{\sigma}^{pk}.
\end{equation} 
\end{prop}

\subsection{}\label{pcour5}
Soient, respectivement,    $\widehat{{\rm{D}}_{A, \sigma }}$,  $\widehat{{\rm{Z}}{\rm{D}}_{A, \sigma }}$, $\widehat{{\rm{Z}}{A}_{A, \sigma }}$ les complétés adiques de ${\rm{D}}_{A, \sigma }$,  ${\rm{Z}}{\rm{D}}_{A, \sigma }$, ${\rm{Z}}{A}_{A, \sigma }$ (\ref{pcour0}) relativement à l'élément central $\partial_{\sigma}^{p} \in {\rm{Z}}{\rm{D}}_{A, \sigma }$. 

 \begin{theo}[\cite{glq2}, Thm. 8.7]\label{frobdiv3}
L'application ${\rm{\Phi}}_{A, \sigma}$   induit  un isomorphisme de $A$-algèbres
\begin{equation}\label{frobdiv31}
\widehat{{\rm{D}}_{A, \sigma }}   \stackrel{\sim} \rightarrow {\rm{End}}_{\widehat{{\rm{Z}}{\rm{D}}_{A, \sigma }}}(\widehat{{\rm{Z}}{A}_{A, \sigma }}).
\end{equation} 
\end{theo}

\subsection{}\label{pcour6}
Pour $q=1$, réduisant modulo $p$,  cet isomorphisme  redonne l'isomorphisme (\ref{sci1}). A un choix de normalisation près (correspondant exactement à la $q$-déformation de la différence entre diviser par $p$ ou par $p!$ dans la construction du Frobenius divisé en caractéristique $p$), pour $A=R[t], f={\rm{Id}}$, l'isomorphisme (\ref{frobdiv31}) se décrit explicitement comme dans (\cite{gl}, \S4).  

\section{Théorie de Hodge non-abélienne $q$-déformée}  
On conserve dans ce \S \, les hypothèses et notations générales du \S\ref{pcour}.
\subsection{}\label{hna1}
Le lemme classique d'algèbre linéaire (\cite{glq}, Lem. 5.6) déjà évoqué en \ref{ov4} montre alors que les anneaux $\widehat{{\rm{D}}_{A, \sigma }}$ et $\widehat{{\rm{Z}}{\rm{D}}_{A, \sigma }}$ sont   équivalents au sens de Morita. On va traduire cette conséquence en termes plus explicites.
\begin{defi}[]\label{sigmader}
Soit $M$ un $A$-module. Une   $\sigma$-dérivation (de niveau 0) ou simplement  {\em{$\sigma$-dérivation}}  de $M$ est une application $R$-linéaire $\partial_{\sigma, M}$ (=: $\partial_{\sigma, M}^{<1>_{0}}$) vérifiant, pour tous $r \in A, m \in M$, l'égalité  (règle de Leibniz $q$-déformée) 
\begin{equation}\label{taylor2}
\partial_{\sigma, M}(rm) = \partial_{\sigma}(r)m + \sigma(r) \partial_{\sigma, M}(m).
\end{equation}
\end{defi}
On a une notion évidente de morphismes entre modules munis de $\sigma$-dérivations.
Rappelons maintenant qu'on dit qu'un endomorphisme $u_G$  d'un groupe abélien $G$ est dit  {\emph{quasi-nilpotent}} si pour tout $g\in G$, il existe $n \in {\Bbb{N}}$ tel que $u_{G}^{n}(g)=0$.

\subsection{}\label{hna2}
Dans le cadre géométrique du  \S\ref{pcour}, l'analogue $q$-déformé de la correspondance d'Ogus-Vologodsky (\ref{cov1}), (\ref{cov2}) est   l'énoncé suivant.
\begin{theo}[\cite{glq2}, Cor. 8.9]\label{qsci}
La catégorie des $A$-modules $M$ munis d'une   $\sigma$-dérivation  quasi-nilpotente $\sigma_M$ est équivalente à la catégorie des $A'$-modules $H$ munis d'un endomorphisme $A$-linéaire quasi-nilpotent $u_H$. 
\end{theo}
L'équivalence est donnée explicitement (comparer avec \cite{glq}, Prop. 5.7 pour la situation en caractéristique $p$) par les deux foncteurs suivants quasi-inverses l'un de l'autre
\begin{equation}
{\Bbb{H}}_{q} : (M, \sigma_M)  \mapsto (H := \{m \in M \, | \,{\rm{\Phi}}_{A, \sigma}(\partial_{\sigma}^{k})(m) = \partial_{\sigma}^{k}(m)\, {\rm{pour \, tout}} \,\,k\in {\Bbb{N}}\}, \partial_{\sigma}^p),
\end{equation}
\begin{equation}
{\Bbb{M}}_{q} :  (H, u_H) \mapsto (M:= A  \,\, \otimes_{A'} \,\, N, \partial_{\sigma, M})
\end{equation}
avec $\partial_{\sigma, M}$ l'unique  $\sigma$-dérivation de $M$ telle que $\partial_{\sigma, M}(1 \otimes h) = t^{p-1}\otimes u_H(h)$ pour tout $h\in H$.
Dans cette équivalence, $(A, \partial_{\sigma})$ (\ref{op5}) correspond à $(A', 0)$.

 \subsection{}\label{hna2}
 Formulons maintenant les conséquences cohomologiques de cette équivalence en termes analogues à ceux de \ref{comp}. 
 Si $M$ est un $A$-module muni d'une $\sigma$-dérivation  $\partial_{\sigma, M}$, on lui associe son complexe de De Rham $q$-déformé ou, s'inspirant de  la terminologie de \cite{scholze},   {\emph{complexe de $q$-De Rham de $M$}}  
 \begin{equation}
q{\text{-}}{\rm{DR}}(M/R) : 0 \rightarrow  M \stackrel{\nabla_{M}} \longrightarrow M \otimes_A  \Omega^{1}_{A/R}\rightarrow 0
\end{equation}
avec $M$ placé en degré 0 et $\nabla_M(m) = \partial_{\sigma, M} (m) \otimes dx$.  Bien que cela ne joue pas de rôle à ce niveau, signalons ici qu'il serait beaucoup plus canonique dans cette définition d'utiliser  le $R$-{\emph{module des différentielles $q$-déformées}} $\Omega^{1}_{A/R, \sigma}$ de \cite{lq2}, Def. 5.3. plutôt que $ \Omega^{1}_{A/R}$ (qui   lui est  seulement non-canoniquement isomorphe).

D'autre part, pour $H$ un $A'$-module muni d'un endomorphisme $u_H$, on peut lui associer son complexe de Higgs 
 \begin{equation}
{\rm{Higgs}}\,(H/R) : 0 \rightarrow  H \stackrel{\theta_H} \longrightarrow H \otimes_{A'}  \Omega^{1}_{A'/R}\rightarrow 0
\end{equation}
avec $H$ placé en degré 0 et $\theta_H (h) = u_H(h) \otimes dx$.
 
\begin{prop}[\cite{glq2}, Cor. 8.10]\label{qcomp}
Si $(M, \partial_{\sigma, M})$ est un $A$-module muni d'une dérivation $q$-déformée quasi-nilpotente et $(H, u_H)$ un $A'$-module muni d'un endomorphisme quasi-nilpotent se correspondant suivant les foncteurs  ${\Bbb{H}}_{q}$ et ${\Bbb{M}}_{q}$, alors le  complexe $q{\text{-}}{\rm{DR}}(M/R) $ est quasi-isomorphe au complexe  ${\rm{Higgs}}\,(H/R)$. 
\end{prop}
 \begin{cor}\label{qcart}
Il existe un isomorphisme (``de Cartier q-déformé'') de $R$-modules
\begin{equation}
C_q :  {\rm{H}}^{i} (q{\text{-}}{\rm{DR}}(A/R)) \stackrel{\sim}   \rightarrow  \Omega^{i} _ {A'/R} 
\end{equation}
pour tout $i$.
\end{cor}
Un suivi des différents morphismes permet de vérifier qu'il s'identifie bien à celui donné dans  \cite{scholze}, Prop. 3.4, (iii) pour l'exemple $R= \co_K$ de \ref{op1}.

\subsection{}\label{} La condition de $q$-divisibilité de \ref{pcour0} garde un sens lorsque $p$ est remplacé par une puissance de $p$ et l'hypothèse de $q$-divisibilité  de $R$ correspondante est   cruciale pour généraliser à ce cadre les principaux résultats ci-dessus (cf. \cite{glq2}).  On notera ici qu'elle n'est, en général, pas vérifiée pour $R= \co_K$ et $q\neq1$ une racine $p^n$-ième ($n>1$) de l'unité comme dans  \ref{op1}. 
 
\section{Questions-Travaux en cours}  

\subsection{}\label{q1}Pour ce qui est du lien avec \cite{bs}, les questions que nous nous posons sont toutes celles motivées par l'espoir suivant, dont les termes seront précisés le moment venu  :\\

{\it{L'équivalence de catégories \ref{qsci} est un corollaire de l'explicitation locale d'une équivalence canonique, compatible (à torsion près en général) au passage à la cohomologie,  entre une catégorie  convenable  de cristaux sur un site $q$-cristallin {\rm{(\cite{bs}, 16.2)}} et une autre de cristaux sur un site prismatique {\rm{(\cite{bs}, 4.1)}}.}}\\

Indépendamment de \cite{bs}, une première étape pourrait  consister à reformuler \ref{qsci} comme un cas particulier d'une équivalence entre des catégories de ${\rm{D}}$-modules $q$-déformés convenables, le modèle ``non $q$-déformé'' étant le point de vue proposé  par Shiho (\cite{shiho}, Thm. 3.1) consistant à  voir la correspondance d'Ogus et Vologodsky comme cas particulier d'un résultat plus général. Pour ce faire, il devrait être utile d'introduire  (suivant les mêmes lignes que celles utilisées pour définir (\ref{opdt}))  un {\emph{anneau d'opérateurs différentiels $q$-tordus de niveau -1}}   (avec $q$ ``générique'') déformant  celui introduit par Shiho (\cite{shiho}, \S2) et intervenant dans sa généralisation de \cite{ov}.

Ensuite, dans une seconde étape, pour faire le lien entre \cite{bs} et nos constructions, l'idée la plus naturelle est de généraliser la classique équivalence entre catégories de cristaux et catégories de  $\cD$-modules et sa compatibilité au passage à la cohomologie au cadre des sites évoqués ci-dessus et des anneaux d'opérateurs différentiels $q$-déformés qui leur correspondent.  

Enfin, il restera, dans une  dernière étape, à définir dans un cadre géométrique non nécessairement ``local''  le foncteur canonique entre cristaux qu'on espère pouvoir s'expliciter comme ``Frobenius divisé'' au niveau des algèbres d'opérateurs différentiels $q$-déformés de niveau 0 et -1. Ce foncteur   devrait simplement  être 
celui induit par  le morphisme image inverse déduit du morphisme de sites décrit dans (\cite{bs}, début de la preuve du Thm. 16.17). 

\subsection{}\label{q2}

Donnons, en conservant  les notations de \ref{op1}  et en supposant que $R$ modulo $(p)_q$ soit $q$-divisible,
quelques indications sur la première étape et sur la définition de l'anneau ${\rm{D}}^{(-1)}_{A, \sigma }$ d'opérateurs différentiels $q$-tordus de niveau -1 (d'autres niveaux négatifs, comme dans \cite{shiho}, sont évidemment possibles).  La définition de ${\rm{D}}^{(-1)}_{A, \sigma }$ suit   celle de ${\rm{D}}^{(0)}_{A, \sigma } $ (\ref{defop}) en remplaçant {\emph{formellement}} partout $A \langle  \xi   \rangle_{q,y}$ par $A \langle  \frac{\xi}{(p)_q}   \rangle_{q^p,y}$ 
et en modifiant en conséquence \ref{multdef}, etc. On montre alors que la donnée d'une structure de ${\rm{D}}^{(-1)}_{A, \sigma }$-module sur un  $A$-module $M$ est équivalente à la donnée d'une {\emph{$\sigma^p$-dérivation de niveau -1}}, i.e. (comparer avec \ref{sigmader}) d'une application $R$-linéaire  $\partial_{\sigma^p, M}^{<1>_{1}} :  M \rightarrow M$ telle que, pour tous $r \in A, m \in M$, on ait.
\begin{equation}\label{}
\partial_{\sigma^p, M}^{<1>_{1}} (rm) = (p)_q\partial_{\sigma^p}(r)m + \sigma^p(r) \partial_{\sigma^p, M}^{<1>_{1}} (m).
\end{equation}
Il est facile de voir qu'il existe un foncteur "image inverse par Frobenius (relatif)", analogue $q$-déformé de (\cite{shiho}, Thm. 3.1), de la catégorie des ${\rm{D}}^{(-1)}_{A', \sigma }$-modules dans celle des ${\rm{D}}^{(0)}_{A, \sigma }$-modules dont nous pensons savoir  démontrer (\cite{glq3}) que c'est une équivalence de catégories sur les objets quasi-nilpotents.

 En particulier, lorsque $q^p=1$, un ${\rm{D}}^{(-1)}_{A', \sigma }$-module $M$ n'est pas autre chose qu'un $A'$-module de Higgs et le  théorème \ref{qsci} serait alors un cas particulier de cette  équivalence de catégories plus générale.    

\subsection{}\label{q3}
Reprenons  les notations de \ref{op1} et supposons de plus   que $R$ soit une  algèbre au-dessus de  ${\Bbb{Z}}_{p}[[q-1]]$  munie d'une structure de $\delta$-anneau (\cite{bs}, Def. 2.1)  telle  que $\delta(q)=0$ (comme pour ${\Bbb{Z}}_{p}[[q-1]])$.  Supposons également $A$   munie d'une structure de $\delta$-$R$-algèbre telle que $\delta(x)=0$,  structure qu'on étendra à  $A[\xi]$ en posant  
\begin{equation}\label{}
\delta(\xi) = \sum_{1 \leq i \leq p-1} \frac{1}{p} \left(\begin{array}{c}p \\i\end{array}\right)x^{p-i}\xi^i.
\end{equation} 
Le premier site qui nous intéresse  dans \ref{q1} est le {\emph{site $q$-cristallin}}\footnote{CAVEAT : Nous empruntons   ici et plus bas, abusivement,    la terminologie de  \cite{bs} mais ignorons dans nos rappels certaines des propriétés additionnelles sur les objets requises  dans {\it{loc. cit.}} si elles ne jouent pas de rôle  dans ce que l'on veut expliquer ici (voir d'ailleurs, à ce sujet,  les commentaires sur leur éventuel caractère provisoire sous la définition 16.2 de \cite{bs}). Les ajustements  précis avec les hypothèses de \cite{bs}, particulièrement ceux nécessitant de prendre en compte  complétions et  topologies (ne serait-ce que dans la définition des sites)  seront donnés  dans \cite{glq3}.}(\cite{bs},  16.2) de $A/(q-1)$ relativement à $(R, (q-1))$.
Le  point crucial dans  la seconde étape espérée dans \ref{q1} est une identification    de $(A \langle  \xi   \rangle_{q,y}, I^{[1]_{q,y}})$ (\ref{op3}) avec  la  $q$-PD-enveloppe (\cite{bs}, Lemma 16.10) de $(A[\xi],(\xi))$.  Précisons le résultat auquel nous parvenons.   Pour définition d'une {\emph{q-PD-paire}} $(B,J)$  (comparer \cite{bs}, Def. 16.2) ne retenons ici (cf. $^{1}$) que la donnée  d'une $\delta$-algèbre $B$ au-dessus de $(R, \delta)$, sans $(p)_p$-torsion,    munie  d'un idéal $J$ tel que $\phi(J) \subset (p)_qB$  
avec $\phi (b) := b^p +p\delta(b)$  pour tout $b \in B$. Si $C$ est une $\delta$-$R$-algèbre et $I$ un idéal quelconque de $C$ (auquel cas, on dira que $(C,I)$ est une  {\emph{$\delta$-paire}}, cf. (\cite{bs}, Def. 3.2)), sa {\emph{$q$-PD-enveloppe}}, notée $(C^{[\,]}, I^{[\,]})$,  est pour nous  ici la $q$-PD-paire universelle  (dont la proposition ci-dessous prouve, pour le cas qui la concerne, 
l'existence et l'unicité à isomorphisme près) pour le prolongement  (unique) à $(C^{[\,]}, I^{[\,]})$ de tout morphisme $(C,I) \rightarrow (B,J)$  d'une $\delta$-paire dans une $q$-PD-paire.   On a alors la  
\begin{prop}[\cite{glq3}]\label{q-envel}
Si $A$ est une $\delta$-$R$-algèbre sans $(p)_q$-torsion, alors $(A \langle  \xi   \rangle_{q,y}, I^{[1]_{q,y}})$ {\rm{(\ref{op3})}} s'identifie à la $q$-PD-enveloppe$^{\rm{1}}$ de la $\delta$-paire $(A[\xi], (\xi))$.
\end{prop}

La démonstration  consiste à se ramener au cas $R={\Bbb{Z}}_{p}[[q-1]]$ et $A=R[x]$ puis, utilisant l'écriture $p$-adique de $n = \sum_{r\geq0} k_rp^r$, de montrer que les $v_n := \xi^{k_0} \prod_{r\geq0} (\delta^{r}([\phi](\xi)))^{k_r+1}$ forment, pour $n \in {\Bbb{N}}$,  une base, comme $A$-module,  de $A \langle  \xi   \rangle_{q,y}$. Ici 
\begin{equation}\label{[phi]}
[\phi] : A \langle  \xi   \rangle_{q,y} \rightarrow A \langle  \xi   \rangle_{q,y} \,\,;\,\, [\phi](\xi^{[n]_{q,y}}) = \sum_{i=n}^{pn}B_{n,i}(q)x^{pn-i}\xi^{[i]_{q,y}}
\end{equation}
tient  compte, par rapport à  $[{\rm{F}}^*]$ (\ref{frodiv1}), de l'usage du Frobenius absolu  dans \cite{bs} plutôt que relatif dans \cite{glq2}
. Si maintenant $(B,J)$ est une $q$-PD-paire, tout morphisme de $\delta$-paires $u : (A[\xi], (\xi)) \rightarrow (B,J)$ s'étend alors uniquement à $A \langle  \xi   \rangle_{q,y}$ par un morphisme d'anneaux  envoyant $v_n \in A \langle  \xi   \rangle_{q,y}$ sur $f^{k_0} \prod_{r\geq0} (\delta^{r}(g))^{k_r+1} \in B$ avec $f:=u(\xi)$ et  $g\in B$ unique  tel que $\phi(f) = (p)_qg$.  \\

Signalons que le cas $q=1$ est celui traité dans  (\cite{bs}, Lem. 2.35) et, pour le lecteur intéressé, l'existence d'un analogue (\cite{prid}, Lem. 1.3), au moins lorsque $q-1 \in R^{\times}$, pour les $\lambda$-anneaux.  

\subsection{}\label{q4}
Conservons les notations de \ref{q3} et supposons que    $(R, (p)_q)$ soit un {\emph{prisme borné}} (\cite{bs}, Def. 3.2) pour pouvoir  réfèrer à \cite{bs}. Notons $A^{(1)}$ le
quotient $A'/(p)_qA'$. Le second site qui nous intéresse dans \ref{q1} est le {\emph{site prismatique}}  de $A^{(1)}$ relativement à $(R, ((p)_q))$.
Comme dans \ref{q3}, le point crucial dans  la seconde étape espérée dans \ref{q1} est  de disposer d'une description adéquate (que nous appliquerons {\it{in fine}} à $A'$ plutôt qu'à $A$) de l'{\emph{enveloppe prismatique}} (\cite{bs}, Cor. 3.14) de $(A[\xi], (\xi))$ relativement à $(R, ((p)_q))$. 
Reprenant les termes de la construction  donnée dans {\it{loc. cit.}},  considérons donc  juste ici la question du prolongement universel  d'un morphisme de  $\delta$-paires (au-dessus de la $\delta$-paire $(R, ((p)_q))$)  $u : (A[\xi], (\xi))  \rightarrow (B, ((p)_q))$ avec $B$ sans $(p)_q$-torsion à   une $\delta$-paire de la forme $(C, ((p)_q))$. Nous montrons  qu'un tel objet universel existe et nous l'appellerons  (cf. $^1$) dans la proposition qui suit {\emph{enveloppe prismatique}} de $(A[\xi], (\xi))$.
\begin{prop}[\cite{glq3}]\label{prism-envel}
Si $A$ est une $\delta$-$R$-algèbre sans $(p)_q$-torsion, alors $(A  \langle  \frac{\xi}{(p)_q}   \rangle_{q^p,y}, ((p)_q))$ {\rm{(\ref{q2})}} s'identifie à l'enveloppe prismatique$^{\rm{1}}$  de  la $\delta$-paire $(A [\xi], (\xi))$.
\end{prop}
En effet, la variante de (\ref{frobeniustordu})  utilisant le Frobenius absolu de $A$ plutôt que relatif comme dans \cite{glq2} fournit une application (dont (\ref{[phi]}) est la variante ``divisée'')
\begin{equation}\label{phi}
 \phi  : A \langle  \xi   \rangle_{q,y} \rightarrow A \langle  \xi   \rangle_{q,y} \,\,;\,\,  \phi (\xi^{[n]_{q,y}}) = \sum_{i=n}^{pn}A_{n,i}(q)x^{pn-i}\xi^{[i]_{q,y}}.
 \end{equation}
L'équation (\ref{defB}) suffit alors à voir (rappelons au passage que $\phi(q)=q^p$) que $A \langle  \frac{\xi}{(p)_q}   \rangle_{q^p,y}$ est bien muni d'un relèvement de Frobenius et, par suite, d'une structure de $\delta$-anneau.
Enfin, le même argument que dans \ref{q-envel} (avec $u(\xi) =  (p)_qg$) donne le prolongement cherché de $u : (A[\xi], (\xi))  \rightarrow (B, (p)_qB)$ à $(A \langle  \frac{\xi}{(p)_q}   \rangle_{q^p,y}, ((p)_q))$.

\subsection{}\label{q5}
Pour terminer, remarquons que les arguments de Shiho (\cite{shiho}) ne   nécessitaient pas d'interprétation de ses ${\rm{D}}^{(-1)}$-modules ({\it{loc. cit}} \S2) quasi-nilpotents en termes de cristaux sur un site mais que, lorsque $q=1$,  le  site prismatique (\cite{bs}, 4.1) en fournit une, qui dans ce cas particulier est juste une variante ``avec $\delta$-structures''   de celle   déjà établie  dans (\cite{oyama}, Def. 1.3.1, \cite{xu}, Def. 7.1)  (l'ajout de $\delta$-structures évitant  précisément les puissances divisées additionnelles sur les anneaux d'opérateurs différentiels considérés des ces articles). Enfin, compte tenu des considérations topologiques délicates à développer sur les sites considérés dans   \ref{q3}-\ref{q4} nous laissons pour ailleurs la discussion  d'une possible approche alternative  à l'équivalence cherchée dans \ref{q1} qui serait l'analogue de   (\cite{oyama}, Thm. 1.4.3) (équivalence de topos).

%


\begin{thebibliography}{99} 



\bibitem{agt} {\sc A. Abbes, M. Gros, T. Tsuji},  The $p$-adic Simpson correspondence.  Annals of Math. Stud., 193, Princeton Univ. Press (2016).

\bibitem{ber} {\sc P. Berthelot}, {\em $\cD$-modules arithmétiques. I. Opérateurs différentiels de niveau fini}. Ann. Sc. de l'\'E.N.S, 4 ième série, t. 29, n. 2 (1996), 185-272.

\bibitem{bms} {\sc B. Bhatt, M. Morrow, P. Scholze}, {\em Integral p-adic Hodge theory}. Publ. Math. Inst. Hautes Études Sci. 128 (2018), 219-397.



 \bibitem{bs} {\sc B. Bhatt,  P. Scholze}, {\em Prisms and prismatic cohomology}. arXiv:1905.08229v2 (2019).


\bibitem{di} {\sc P. Deligne, L. Illusie}, {\em Relèvements modulo $p^2$ et décomposition du complexe de de Rham}.  Inv. Math.  {\bf 89} (1987),  247-270.

\bibitem{glq} {\sc M. Gros, , B. Le Stum, A. Quir\'os}, {\em A Simpson correspondence in positive characteristic}. Pub. RIMS Kyoto Univ.  {\bf 46} (2010), 1-35. 

\bibitem{gl} {\sc M. Gros, , B. Le Stum}, {\em Une neutralisation explicite de l'algèbre de Weyl quantique complétée}.  Comm. Algebra  {\bf 42}  no. 5, (2014), 2163-2170. 



\bibitem{glq2} {\sc M. Gros, , B. Le Stum, A. Quir\'os}, {\em  Twisted divided powers and applications}. J. Number Theory (2019), (https://doi.org/10.1016/j-jnt.2019.02.009).

\bibitem{glq3} {\sc M. Gros, , B. Le Stum, A. Quir\'os}, {\em  En préparation}.  

\bibitem{katz} {\sc N. Katz}, {\em Nilpotent connections and the monodromy theorem : applications of a result of Turrittin}. Publ. Math. IHÉS {\bf 39} (1970), 175-232. 


\bibitem{lq2} {\sc  B. Le Stum, A. Quir\'os}, {\em Twisted calculus}. Communications in Algebra, 46:12 (2018), 5290-5319. (https://doi.org/10.1080/00927872.2018.1464168). 

\bibitem{lq} {\sc   B. Le Stum, A. Quir\'os}, {\em Formal confluence of quantum differential operators}. Pacific Journal of Mathematics {\bf 292(2)}, (2018), 427-478.  (https://doi.org/10.2140/pjm.2018.292.427).
 
\bibitem{ov} {\sc A. Ogus, V. Vologodsky},  {\em Non abelian  Hodge theory in characteristic $p$}. Publ. Math. IHÉS {\bf 101} (2007), 1-138. 

\bibitem{oyama} {\sc H. Oyama},  {\em PD Higgs crystals and Higgs cohomology in characteristic $p$}.  J. Algebraic Geom. {\bf 26} (2017), 735-802.
 
 \bibitem{prid} {\sc J. P. Pridham},  {\em On $q$-De Rham cohomology via $\Lambda$-rings}. Math. Annalen  (https://doi.org/10.1007/s00208-019-01806-7). 

\bibitem{scholze} {\sc P. Scholze},   {\em Canonical $q$-deformations in arithmetic geometry.} Ann. Fac. Sci. Toulouse Math. (6) 26 (2017), no. 5, 1163-1192.  
 
 
 


\bibitem{shiho} {\sc A. Shiho}, {\em Notes on generalizations of local Ogus-Vologodsky correspondence.}  J. Math. Sci. Univ. Tokyo  {\bf 22} (2015), 793-875. 

\bibitem{tsuji} {\sc T. Tsuji}, {\em Exposé du 12 Mai 2017 au Simons Symposium on $p$-adic Hodge theory  et articles  en préparation dont  un en collaboration avec Matthew Morrow.}  


\bibitem{xu} {\sc D. Xu}, {\em Lifting the Cartier transform of Ogus-Vologodsky modulo $p^n$}.  Mémoires de la Société Mathématique de France, t. {\bf{163}} (2019).


\end{thebibliography}
\end{document}